\documentclass[10pt,reqno]{amsart}
\usepackage{amssymb}
\usepackage{amsxtra}
\usepackage[mathscr]{eucal}

\theoremstyle{plain}

\newtheorem{Prop}{Proposition}
\theoremstyle{definition}
\newtheorem{Def}{Definition}
\theoremstyle{remark}

\begin{document}

\title[David M.~Bradley]{A Signed Analog of Euler's Reduction Formula
for the Double Zeta Function}

\date{\today}

\author{David~M. Bradley}
\address{Department of Mathematics \& Statistics\\
         University of Maine\\
         5752 Neville Hall
         Orono, Maine 04469-5752\\
         U.S.A.}
\email[]{bradley@math.umaine.edu, dbradley@member.ams.org}

\subjclass{Primary: 11M41; Secondary: 11M06}

\keywords{Euler sums, multiple harmonic series, double zeta values.}

\begin{abstract}
The double zeta function is a function of two arguments defined by a
double Dirichlet series, and was first studied by Euler in response
to a letter from Goldbach in 1742.  By calculating many examples,
Euler inferred a closed form evaluation of the double zeta function
in terms of values of the Riemann zeta function, in the case when
the two arguments are positive integers with opposite parity.  Here,
we consider a signed analog of Euler's evaluation: namely a
reduction formula for the signed double zeta function that reduces
to Euler's evaluation when the signs are specialized to 1.  This
formula was first stated in a 1997 paper by Borwein, Bradley and
Broadhurst and was subsequently proved by Flajolet and Salvy using
contour integration.  The purpose here is to give an elementary
proof based on a partial fraction identity.
\end{abstract}

\maketitle

%\tableofcontents
\interdisplaylinepenalty=500

\section{Introduction}\label{sect:Intro}
%The Riemann zeta function is defined for $\Re(s)>1$ by
%\begin{equation}\label{Riemann}
%   \zeta(s) := \sum_{n=1}^\infty \frac{1}{n^s}.
%\end{equation}
The double zeta function is defined by
\begin{equation}\label{doublezeta}
   \zeta(s,t)
   := \sum_{n=1}^\infty \frac{1}{n^s}\sum_{k=1}^{n-1}\frac{1}{k^t},
      \qquad \Re(s)>1, \qquad \Re(s+t)>2.
\end{equation}
The problem of evaluating sums of the form~\eqref{doublezeta} with
integers $s>1$ and $t>0$ seems to have been first proposed in a
letter from Goldbach to Euler~\cite{LE2} in 1742. (See
also~\cite{LE,Goldbach} and~\cite[p.\ 253]{Berndt1}.)  Calculating
several examples led Euler to infer a closed form evaluation of the
double zeta function in terms of values of the Riemann zeta function
\[
   \zeta(s) := \sum_{n=1}^\infty \frac{1}{n^s},\qquad \Re(s)>1,
\]
in the case when $s-1$ and $t-1$ are positive integers with opposite
parity. Euler's evaluation can be expressed as follows. Let $s>1$
and $t>1$ be integers with opposite parity (i.e.\ $s+t$ is odd) and
let $2m=\max(s,t)$. Then
\begin{multline}\label{EulerRed}
   \zeta(s,t) = \frac12\big((1+(-1)^s\big)\zeta(s)\zeta(t)+\frac12\bigg[(-1)^s
                \binom{s+t}{s}-1\bigg]\zeta(s+t)\\
               +(-1)^{s+1} \sum_{k=1}^{m}\bigg[\binom{s+t-2k-1}{t-1}+
               \binom{s+t-2k-1}{s-1}\bigg]\zeta(2k)\zeta(s+t-2k).
\end{multline}
The formula~\eqref{EulerRed} is also valid when $t=1$ and $s$ is
even, but that case is subsumed by another formula of Euler,
namely
\begin{equation}\label{zs1}
   \zeta(s,1) = \frac{1}{2} s\zeta(s+1) - \frac12 \sum_{k=2}^{s-1}
   \zeta(k)\zeta(s+1-k),
\end{equation}
which is valid for all integers $s>1$.  In~\cite{BBB}, Borwein,
Bradley and Broadhurst considered the more general Euler sum
\begin{equation}\label{EulerSum}
   \zeta(s_1,s_2,\dots,s_k;\sigma_1,\sigma_2,\dots,\sigma_k) :=
   \sum_{n_1>n_2>\cdots>n_k>0}\; \prod_{j=1}^k \sigma_j^{n_j} n_j^{-s_j}
\end{equation}
with each $\sigma_j\in\{-1,1\}$.   Among the many other results
for~\eqref{EulerSum} listed therein is an explicit formula for the
case $k=2$ that reduces to~\eqref{EulerRed} when
$\sigma_1=\sigma_2=1$.  We restate this result as follows:
\begin{Prop}\label{prop:BBB}
  Let $\sigma,\tau\in\{-1,1\}$, and
  let $s$ and $t$ be positive integers such that $s+t$ is odd,
  $s>(1+\sigma)/2$, and $t>(1+\tau)/2$.  Then
  \begin{multline}\label{BBBeq75}
    \zeta(s,t;\sigma,\tau)
    = \tfrac12\big(1+(-1)^s\big)\zeta(s;\sigma)\zeta(t;\tau)
    -\tfrac12\zeta(s+t;\sigma\tau)\\
    +(-1)^{t}\sum_{0\le k\le t/2}
    \binom{s+t-2k-1}{s-1}\zeta(2k;\sigma\tau)\zeta(s+t-2k;\sigma)\\
    +(-1)^{t}\sum_{0\le k\le s/2}\binom{s+t-2k-1}{t-1}\zeta(2k;\sigma\tau)\zeta(s+t-2k;\tau).
  \end{multline}
\end{Prop}
In Proposition 1, it is understood that $\zeta(0;\sigma\tau)=-1/2$
in accordance with the analytic continuation of $s\mapsto
\zeta(s;\sigma\tau)$.  The restriction $t>(1+\tau)/2$ can be removed
if in~\eqref{BBBeq75} we interpret $\zeta(1;1)=0$ wherever it
occurs.  That is, if $\sigma\in\{-1,1\}$ and $s$ is an even positive
integer, then
\begin{equation}\label{BBBeq75t=1}
   \zeta(s,1;\sigma,1) = \frac12
   (s-1)\zeta(s+1;\sigma)+\frac12\zeta(s+1)-\sum_{k=1}^{(s/2)-1}
   \zeta(2k;\sigma)\zeta(s+1-2k).
\end{equation}
Note that by~\eqref{zs1} we know that the case $\sigma=1$
of~\eqref{BBBeq75t=1} can be extended to \emph{all} integers $s>1$,
not just even $s$.

Using contour integration, Flajolet and Salvy~\cite{FlajSalv} proved
an equivalent version of Proposition~\ref{prop:BBB}.  Our intention
here is to give an elementary proof based on the partial fraction
decomposition
\begin{equation}\label{parfrac}
   \frac{1}{x^sy^t} = \sum_{a=0}^{s-1} \binom{a+t-1}{t-1}
   \frac{1}{x^{s-a}(x+y)^{t+a}} +
   \sum_{a=0}^{t-1}\binom{a+s-1}{s-1}\frac{1}{(x+y)^{s+a}y^{t-a}},
\end{equation}
which is valid for positive integers $s$ and $t$ and non-zero real
numbers $x$ and $y$ such that $x+y\ne 0$.  As in~\cite{DBqDecomp},
we note that~\eqref{parfrac} is readily proved by applying the
partial differential operator
\[
   \frac{1}{(r-1)!}\bigg(-\frac{\partial}{\partial x}\bigg)^{r-1}
   \frac{1}{(s-1)!}\bigg(-\frac{\partial}{\partial y}\bigg)^{s-1}
\]
to both sides of the identity
\[
   \frac{1}{xy}
 = \frac{1}{x+y}\bigg(\frac{1}{x}+\frac{1}{y}\bigg).
\]

\section{Proof of Proposition~\ref{prop:BBB}}
\begin{Def} Let $N$ be a positive integer and let $s,t,\sigma,\tau$
be complex numbers.  Define
\[
   \zeta_N(s,t;\sigma,\tau)
   := \sum_{n=1}^N \sum_{k=1}^{n-1}\frac{\sigma^n \tau^k}{n^sk^t}
   % = \sum_{n=1}^N\sum_{k=1}^{n-1}\frac{\sigma^n\tau^{n-k}}{n^s(n-k)^t}
   =(-1)^t\sum_{n=1}^N\sum_{k=1}^{n-1}\frac{\sigma^n\tau^{n-k}}{n^s(k-n)^t}
   \quad\text{and}\quad
   \zeta_N(s;\sigma) = \sum_{n=1}^N \frac{\sigma^n}{n^s}.
\]
\end{Def}
Now suppose that $s$ and $t$ are positive integers.
In~\eqref{parfrac} let $x=n$, $y=k-n$, multiply through by $(-1)^t
\sigma^n \tau^{n-k}$ and sum over all positive integers $n$ and $k$
satisfying $N>n>k>0$.  We find that
\begin{align*}
   &(-1)^t\zeta_N(s,t;\sigma,\tau)\\ &= \sum_{a=0}^{s-1}\binom{a+t-1}{t-1}\sum_{n=1}^N\sum_{k=1}^{n-1}
   \frac{\sigma^n\tau^{n-k}}{n^{s-a}k^{t+a}}+\sum_{a=0}^{t-1}\binom{a+s-1}{s-1}
   \sum_{n=1}^N\sum_{k=1}^{n-1}\frac{\sigma^n\tau^{n-k}}{k^{s+a}(k-n)^{t-a}}\\
   &= \sum_{a=0}^{s-1}\binom{a+t-1}{t-1}\bigg\{\sum_{n=1}^N\sum_{k=1}^N
   \frac{\sigma^n\tau^{n-k}}{n^{s-a}k^{t+a}}-\sum_{k=1}^N\frac{\sigma^k}{k^{s+t}}
   -\sum_{k=1}^N\sum_{n=1}^{k-1}\frac{\sigma^n\tau^{n-k}}{n^{s-a}k^{t+a}}\bigg\}\\
   &+(-1)^t\sum_{a=0}^{t-1}\binom{a+s-1}{s-1}(-1)^a\sum_{n=1}^N\sum_{k=1}^{n-1}\frac{\sigma^n\tau^{n-k}}{k^{s+a}(n-k)^{t-a}}
   \\
   &=\sum_{a=0}^{s-1}\binom{a+t-1}{t-1}
   \big[\zeta_N(s-a;\sigma\tau)\zeta_N(t+a;1/\tau)-\zeta_N(t+a,s-a;1/\tau,\sigma\tau)\big]\\
   &-\binom{s+t-1}{s-1}\zeta_N(s+t;\sigma)+(-1)^t\sum_{a=0}^{t-1}\binom{a+s-1}{s-1}(-1)^a\sum_{k=1}^{N-1}\sum_{m=1}^{N-k}
   \frac{\sigma^{m+k}\tau^m}{k^{s+a}m^{t-a}}.
\end{align*}
It follows that
\begin{align}\label{symmetric}
 &(-1)^{t}\zeta_{N}(s,t;\sigma,\tau)+(-1)^s\zeta_{N}(t,s;\tau,\sigma)\nonumber\\
 &\qquad=\sum_{a=0}^{s-1}\binom{a+t-1}{t-1}\zeta_{N}(s-a;\sigma\tau)\zeta_{N}(t+a;1/\tau)
 +\sum_{a=0}^{t-1}\binom{a+s-1}{s-1}\zeta_{N}(t-a;\sigma\tau)\zeta_{N}(s+a;1/\sigma)\nonumber\\
 &\qquad-\bigg[\sum_{a=0}^{s-1}\binom{a+t-1}{t-1}\zeta_{N}(t+a,s-a;1/\tau,\sigma\tau)
 +\sum_{a=0}^{t-1}\binom{a+s-1}{s-1}\zeta_{N}(s+a,t-a;1/\sigma,\sigma\tau)\bigg]\nonumber\\
 &\qquad-\binom{s+t-1}{s-1}\zeta_{N}(s+t;\sigma)-\binom{s+t-1}{t-1}\zeta_{N}(s+t;\tau)\nonumber\\
 &\qquad+(-1)^{t}\sum_{a=0}^{t-1}\binom{a+s-1}{s-1}(-1)^a
 \sum_{k=1}^{N-1}\,\sum_{m=1}^{N-k}\frac{\sigma^{m+k}\tau^m}{k^{s+a}m^{t-a}}\nonumber\\
 &\qquad+(-1)^s\sum_{a=0}^{s-1}\binom{a+t-1}{t-1}(-1)^a
 \sum_{k=1}^{N-1}\,\sum_{m=1}^{N-k}\frac{\tau^{m+k}\sigma^m}{k^{t+a}m^{s-a}}.
\end{align}
By~\eqref{parfrac} again,
\begin{align*}
  &\zeta_N(s;1/\sigma)\zeta_N(t;1/\tau) = \sum_{x=1}^N\sum_{y=1}^N\frac{\sigma^{-x}\tau^{-y}}{x^sy^t}\nonumber\\
%  =\sum_{n=1}^{N}\sum_{x=1}^{n-1}\frac{\sigma^{-x}\tau^{-(n-x)}}{x^s(n-x)^t}
%  +\sum_{n=N+1}^{2N}\,\sum_{x=n-N}^N\frac{\sigma^{-x}\tau^{-(n-x)}}{x^s(n-x)^t}\nonumber\\
  &=
  \sum_{a=0}^{s-1}\binom{a+t-1}{t-1}\sum_{x=1}^{N}\sum_{y=1}^N\frac{\sigma^{-x}\tau^{-y}}{x^{s-a}(x+y)^{t+a}}
  +\sum_{a=0}^{t-1}\binom{a+s-1}{s-1}\sum_{x=1}^{N}\sum_{y=1}^N\frac{\sigma^{-x}\tau^{-y}}{(x+y)^{s+a}y^{t-a}}\nonumber\\
  &=\sum_{a=0}^{s-1}\binom{a+t-1}{t-1}\bigg[\sum_{n=1}^{N}\sum_{x=1}^{n-1}\frac{\sigma^{-x}\tau^{x-n}}{x^{s-a}n^{t+a}}
  +\sum_{n=N+1}^{2N}\,\sum_{x=n-N}^N\frac{\sigma^{-x}\tau^{x-n}}{x^{s-a}n^{t+a}}\bigg]\nonumber\\
  &+\sum_{a=0}^{t-1}\binom{a+s-1}{s-1}\bigg[\sum_{n=1}^{N}\sum_{y=1}^{n-1}\frac{\sigma^{y-n}\tau^{-y}}{n^{s+a}y^{t-a}}
  +\sum_{n=N+1}^{2N}\,\sum_{y=n-N}^N\frac{\sigma^{y-n}\tau^{-y}}{n^{s+a}y^{t-a}}\bigg]\nonumber\\
  &=\sum_{a=0}^{s-1}\binom{a+t-1}{t-1}\bigg[\zeta_{N}(t+a,s-a;1/\tau,\tau/\sigma)
  +\sum_{n=N+1}^{2N}\,\sum_{x=n-N}^N\frac{\sigma^{-x}\tau^{x-n}}{x^{s-a}n^{t+a}}\bigg]\nonumber\\
  &+\sum_{a=0}^{t-1}\binom{a+s-1}{s-1}\bigg[\zeta_{N}(s+a,t-a;1/\sigma,\sigma/\tau)
  +\sum_{n=N+1}^{2N}\,\sum_{y=n-N}^N\frac{\sigma^{y-n}\tau^{-y}}{n^{s+a}y^{t-a}}\bigg].
\end{align*}
Rearranging this yields
\begin{multline}\label{shuffleproduct}
   \sum_{a=0}^{s-1}\binom{a+t-1}{t-1}\zeta_{N}(t+a,s-a;1/\tau,\tau/\sigma)
   +\sum_{a=0}^{t-1}\binom{a+s-1}{s-1}\zeta_{N}(s+a,t-a;1/\sigma,\sigma/\tau)\\
   =\zeta_N(s;1/\sigma)\zeta_N(t;1/\tau)
   -\sum_{a=0}^{s-1}\binom{a+t-1}{t-1}\sum_{n=N+1}^{2N}\,\sum_{x=n-N}^N\frac{\sigma^{-x}\tau^{x-n}}{x^{s-a}n^{t+a}}\\
   -\sum_{a=0}^{t-1}\binom{a+s-1}{s-1}\sum_{n=N+1}^{2N}\,\sum_{y=n-N}^N\frac{\sigma^{y-n}\tau^{-y}}{n^{s+a}y^{t-a}}.
\end{multline}
 Henceforth assume that $\sigma\in\{-1,1\}$
and $\tau\in\{-1,1\}$. Then $\sigma=1/\sigma$, $\tau=1/\tau$ and
from~\eqref{symmetric} and~\eqref{shuffleproduct} we infer that
\begin{align}\label{stillsymmetric}
  &  (-1)^t\zeta_{N}(s,t;\sigma,\tau)+(-1)^s\zeta_{N}(t,s;\tau,\sigma)\nonumber\\
  &=\sum_{a=0}^{s-1}\binom{a+t-1}{t-1}\zeta_{N}(s-a;\sigma\tau)\zeta_{N}(t+a;\tau)
 +\sum_{a=0}^{t-1}\binom{a+s-1}{s-1}\zeta_{N}(t-a;\sigma\tau)\zeta_{N}(s+a;\sigma)\nonumber\\
 &- \zeta_N(s;\sigma)\zeta_N(t;\tau)
 -\binom{s+t-1}{s-1}\zeta_{N}(s+t;\sigma)-\binom{s+t-1}{t-1}\zeta_{N}(s+t;\tau)\nonumber\\
 &+(-1)^{t}\sum_{a=0}^{t-1}\binom{a+s-1}{s-1}(-1)^a
 \sum_{k=1}^{N-1}\,\sum_{m=1}^{N-k}\frac{\sigma^{m+k}\tau^m}{k^{s+a}m^{t-a}}
 +(-1)^s\sum_{a=0}^{s-1}\binom{a+t-1}{t-1}(-1)^a
 \sum_{k=1}^{N-1}\,\sum_{m=1}^{N-k}\frac{\tau^{m+k}\sigma^m}{k^{t+a}m^{s-a}}\nonumber\\
 &+\sum_{a=0}^{s-1}\binom{a+t-1}{t-1}\sum_{n=N+1}^{2N}\,\sum_{x=n-N}^N\frac{\sigma^{x}\tau^{x+n}}{x^{s-a}n^{t+a}}
 +\sum_{a=0}^{t-1}\binom{a+s-1}{s-1}\sum_{n=N+1}^{2N}\,\sum_{y=n-N}^N\frac{\sigma^{y+n}\tau^{y}}{n^{s+a}y^{t-a}}.
\end{align}
But if $N$ is any positive integer, then
\begin{align*}
   \zeta_N(s;\sigma)\zeta_N(t;\tau) &=
   \sum_{x=1}^N\sum_{y=1}^N\frac{\sigma^x\tau^y}{x^sy^t} =
   \sum_{x=1}^N\sum_{y=1}^{x-1}\frac{\sigma^x\tau^y}{x^sy^t}+\sum_{x=y=1}^N\frac{\sigma^x\tau^y}{x^sy^t}
   +\sum_{y=1}^N\sum_{x=1}^{y-1}\frac{\sigma^x\tau^y}{x^sy^t}\\
   &=\zeta_N(s,t;\sigma,\tau)+\zeta_N(s+t;\sigma\tau)+\zeta_N(t,s;\tau,\sigma),
\end{align*}
whence
\[
   \zeta_{N}(t,s;\tau,\sigma) =
   \zeta_{N}(s;\sigma)\zeta_{N}(t;\tau)-\zeta_{N}(s+t;\sigma\tau)-\zeta_{N}(s,t;\sigma,\tau).
\]
If we use this in~\eqref{stillsymmetric} and multiply through by
$(-1)^t$ there comes
\begin{align*}
  &\zeta_{N}(s,t;\sigma,\tau)+(-1)^{s+t}\big[\zeta_{N}(s;\sigma)\zeta_{N}(t;\tau)-\zeta_{N}(s+t;\sigma\tau)
  -\zeta_{N}(s,t;\sigma,\tau)\big]\\
  &=(-1)^t\bigg[\sum_{a=0}^{s-1}\binom{a+t-1}{t-1}\zeta_{N}(s-a;\sigma\tau)\zeta_{N}(t+a;\tau)
  +\sum_{a=0}^{t-1}\binom{a+s-1}{s-1}\zeta_{N}(t-a;\sigma\tau)\zeta_{N}(s+a;\sigma)\bigg]\\
  &+(-1)^{t+1}\bigg[\zeta_N(s;\sigma)\zeta_N(t;\tau)+\binom{s+t-1}{s-1}\zeta_{N}(s+t;\sigma)
  +\binom{s+t-1}{t-1}\zeta_{N}(s+t;\tau)\bigg]\\
  &+\sum_{a=0}^{t-1}\binom{a+s-1}{s-1}(-1)^a\sum_{k=1}^{N-1}\,\sum_{m=1}^{N-k}\frac{\sigma^{m+k}\tau^m}{k^{s+a}m^{t-a}}
  +(-1)^{s+t}\sum_{a=0}^{s-1}\binom{a+t-1}{t-1}(-1)^a\sum_{k=1}^{N-1}\,\sum_{m=1}^{N-k}\frac{\tau^{m+k}\sigma^m}{k^{t+a}m^{s-a}}\\
  &+(-1)^t\bigg[\sum_{a=0}^{s-1}\binom{a+t-1}{t-1}\sum_{n=N+1}^{2N}\,\sum_{x=n-N}^N\frac{\sigma^{x}\tau^{x+n}}{x^{s-a}n^{t+a}}
 +\sum_{a=0}^{t-1}\binom{a+s-1}{s-1}\sum_{n=N+1}^{2N}\,\sum_{y=n-N}^N\frac{\sigma^{y+n}\tau^{y}}{n^{s+a}y^{t-a}}\bigg].
\end{align*}
Now assume also that $s+t$ is odd.
%\begin{align*}
%  &2\zeta_{N}(s,t;\sigma,\tau)
%  =\zeta_{N}(s;\sigma)\zeta_{N}(t;\tau)-\zeta_{N}(s+t;\sigma\tau)\\
%  &+(-1)^t\bigg[\sum_{a=0}^{s-1}\binom{a+t-1}{t-1}\zeta_{N}(s-a;\sigma\tau)\zeta_{N}(t+a;\tau)
%  +\sum_{a=0}^{t-1}\binom{a+s-1}{s-1}\zeta_{N}(t-a;\sigma\tau)\zeta_{N}(s+a;\sigma)\bigg]\\
%  &+(-1)^{t+1}\bigg[\zeta_N(s;\sigma)\zeta_N(t;\tau)+\binom{s+t-1}{s-1}\zeta_{N}(s+t;\sigma)
%  +\binom{s+t-1}{t-1}\zeta_{N}(s+t;\tau)\bigg]\\
%  &+\sum_{a=0}^{t-1}\binom{a+s-1}{s-1}(-1)^a\sum_{k=1}^{N-1}\,\sum_{m=1}^{N-k}\frac{\sigma^{m+k}\tau^m}{k^{s+a}m^{t-a}}
%  -\sum_{a=0}^{s-1}\binom{a+t-1}{t-1}(-1)^a\sum_{k=1}^{N-1}\,\sum_{m=1}^{N-k}\frac{\tau^{m+k}\sigma^m}{k^{t+a}m^{s-a}}\\
%
%&+(-1)^t\bigg[\sum_{a=0}^{s-1}\binom{a+t-1}{t-1}\sum_{n=N+1}^{2N}\,\sum_{x=n-N}^N\frac{\sigma^{-x}\tau^{x-n}}{x^{s-a}n^{t+a}}
% +\sum_{a=0}^{t-1}\binom{a+s-1}{s-1}\sum_{n=N+1}^{2N}\,\sum_{y=n-N}^N\frac{\sigma^{y-n}\tau^{-y}}{n^{s+a}y^{t-a}}\bigg].
%\end{align*}
Writing $(-1)^s$ for $(-1)^{t+1}$ and re-indexing the sums yields
\begin{align}\label{finite}
  2\zeta_{N}(s,t;\sigma,\tau)
  &=
  \big(1+(-1)^s\big)\zeta_{N}(s;\sigma)\zeta_{N}(t;\tau)-\zeta_{N}(s+t;\sigma\tau)\nonumber\\
  &+(-1)^t\sum_{b=1}^s\binom{s+t-b-1}{t-1}\bigg[\zeta_{N}(b;\sigma\tau)\zeta_{N}(s+t-b;\tau)\nonumber\\
  &\qquad\qquad\qquad\qquad\qquad
  +(-1)^b\sum_{k=1}^{N-1}\,\sum_{m=1}^{N-k}\frac{\tau^{m+k}\sigma^m}{k^{s+t-b}m^{b}}
  +\sum_{n=N+1}^{2N}\,\sum_{x=n-N}^N \frac{\sigma^{x}\tau^{x+n}}{x^bn^{s+t-b}}\bigg]\nonumber\\
  &+(-1)^t\sum_{b=1}^t\binom{s+t-b-1}{s-1}\bigg[\zeta_{N}(b;\sigma\tau)\zeta_{N}(s+t-b;\sigma)\nonumber\\
  &\qquad\qquad\qquad\qquad\qquad
  +(-1)^b\sum_{k=1}^{N-1}\,\sum_{m=1}^{N-k}\frac{\sigma^{m+k}\tau^m}{k^{s+t-b}m^{b}}
  +\sum_{n=N+1}^{2N}\,\sum_{y=n-N}^N \frac{\sigma^{y+n}\tau^{y}}{n^{s+t-b}y^b}\bigg]\nonumber\\
  &-(-1)^t\binom{s+t-1}{s-1}\zeta_{N}(s+t;\sigma)-(-1)^t\binom{s+t-1}{t-1}\zeta_{N}(s+t;\tau).
\end{align}
Now suppose that $b$ is odd and $1\le b\le t$.  Then
\begin{align*}
  &\bigg|\zeta_{N}(b;\sigma\tau)\zeta_{N}(s+t-b;\sigma)
  +(-1)^b\sum_{k=1}^{N-1}\,\sum_{m=1}^{N-k}\frac{\sigma^{m+k}\tau^m}{k^{s+t-b}m^b}
  +\sum_{n=N+1}^{2N}\,\sum_{y=n-N}^N\frac{\sigma^{y+n}\tau^{y}}{n^{s+t-b}y^b}\bigg|\\
  &\le
  \bigg|\sum_{k=1}^N\frac{\sigma^k}{k^{s+t-b}}\bigg[\sum_{m=1}^N\frac{(\sigma\tau)^m}{m^b}-
  \sum_{m=1}^{N-k}\frac{(\sigma\tau)^m}{m^b}\bigg]\bigg|
  +\bigg|\sum_{y=1}^N\frac{(\sigma\tau)^y}{y^b}\sum_{n=N+1}^{N+y}
  \frac{\sigma^{n}}{n^{s+t-b}}\bigg|\\
    &\le    \sum_{k=1}^{N}\frac{1}{k^{s+t-b}}\sum_{m=N-k+1}^{N}\frac{1}{m}
    +\sum_{y=1}^N \frac{1}{y}\sum_{n=N+1}^{N+y}\frac{1}{n^{s+t-b}}\\
    &\le \sum_{k=1}^N\frac{1}{k^{s+t-b}}\cdot\frac{k}{N-k+1}
    +\sum_{y=1}^N\frac1y\sum_{n=N+1}^\infty\frac{1}{n^{s+t-b-1}(n-1)}\\
    &=\sum_{k=1}^N\frac{1}{k^{s+t-b-1}(N-k+1)}
    +\sum_{k=1}^N\frac1k\sum_{n=N}^\infty\frac{1}{n(n+1)^{s+t-b-1}}.
\end{align*}
If $t$ is also odd, then (recalling that $s$ and $t$ are positive
integers with opposite parity) $s$ is even and $s+t-b-1\ge s-1 \ge
1$. On the other hand, if $t$ is even, then $b\le t-1$ since $b$ is
odd, and therefore $s+t-b-1 \ge s\ge 1$. In either case,
\begin{multline*}
  \bigg|\zeta_{N}(b;\sigma\tau)\zeta_{N}(s+t-b;\sigma)
  +(-1)^b\sum_{k=1}^{N-1}\,\sum_{m=1}^{N-k}\frac{\sigma^{m+k}\tau^m}{k^{s+t-b}m^b}
  +\sum_{n=N+1}^{2N}\,\sum_{y=n-N}^N\frac{\sigma^{y+n}\tau^{y}}{n^{s+t-b}y^b}\bigg|\\
  \le \sum_{k=1}^{N} \frac{1}{k(N-k+1)}+\sum_{k=1}^N\frac1k\sum_{n=N}^\infty \frac{1}{n(n+1)} =
  \frac{1}{N+1}\sum_{k=1}^{N}\bigg(\frac{1}{k}+\frac{1}{N-k+1}\bigg)+\frac{1}{N}\sum_{k=1}^N\frac1k
  \le \frac{3}{N}\sum_{k=1}^{N}\frac{1}k \to 0
\end{multline*}
as $N\to\infty$.

Next, suppose that $1<b\le t$ and additionally $s>(1+\sigma)/2$.
Then
\[
   \lim_{N\to\infty}\sum_{k=1}^{N-1}\,\sum_{m=1}^{N-k}\frac{\sigma^{m+k}\tau^m}{k^{s+t-b}m^{b}}
   = \sum_{k=1}^\infty \,\sum_{m=1}^\infty
   \frac{\sigma^{m+k}\tau^m}{k^{s+t-b}m^b}
   = \zeta(b;\sigma\tau)\zeta(s+t-b;\sigma)
\]
and
\[
  \bigg|\sum_{n=N+1}^{2N}\,\sum_{y=n-N}^N\frac{\sigma^{y+n}\tau^y}{n^{s+t-b}y^b}\bigg|
  \le \sum_{n=N+1}^{2N}\,\sum_{y=n-N}^N \frac{1}{ny^2}
  =\sum_{y=1}^N\frac{1}{y^2}\sum_{n=N+1}^{N+y}\frac1n
  \le \sum_{y=1}^N\frac{1}{y^2}\cdot\frac{y}{N+1}
  \to 0
\]
as $N\to\infty$.  Thus, under the hypotheses of
Proposition~\ref{prop:BBB}, letting $N\to\infty$ in~\eqref{finite}
yields
\begin{multline*}
   \zeta(s,t;\sigma,\tau) =
   \tfrac12\big(1+(-1)^s\big)\zeta(s;\sigma)\zeta(t;\tau)-\tfrac12\zeta(s+t;\sigma\tau)\\
   +(-1)^t\sum_{1\le k\le s/2}
   \binom{s+t-2k-1}{t-1}\zeta(2k;\sigma\tau)\zeta(s+t-2k;\tau)\\
   +(-1)^t\sum_{1\le k\le t/2}\binom{s+t-2k-1}{s-1}\zeta(2k;\sigma\tau)\zeta(s+t-2k;\sigma)\\
   -\tfrac12(-1)^t\bigg[\binom{s+t-1}{s-1}\zeta(s+t;\sigma)+\binom{s+t-1}{t-1}\zeta(s+t;\tau)\bigg].
\end{multline*}
Since $\zeta(0;\sigma\tau)=-1/2$ by definition, the final line can
be absorbed into the two sums above by permitting $k=0$.

\section{The Case $t=\tau=1$: Proof of Equation~\eqref{BBBeq75t=1}}
Putting $t=\tau=1$ in~\eqref{finite} and noting that then $s$ must
be even yields
\begin{align*}
  2\zeta_N(s,1;\sigma,1) &=
  (s-1)\zeta_N(s+1;\sigma)+\zeta_N(s+1;1)
  +\bigg(\zeta_N(s;\sigma)\zeta_N(1;1)-\sum_{k=1}^{N-1}\sum_{m=1}^{N-k}\frac{\sigma^m}{km^s}\bigg)\\
  &+\bigg(\sum_{k=1}^{N-1}\sum_{m=1}^{N-k}\frac{\sigma^{m+k}}{k^sm}-\zeta_N(s;\sigma)\zeta_N(1;\sigma)\bigg)
   -\bigg(\sum_{n=N+1}^{2N}\,\sum_{x=n-N}^N\frac{\sigma^x}{x^sn}
         +\sum_{n=N+1}^{2N}\,\sum_{y=n-N}^N\frac{\sigma^{y+n}}{n^sy}\bigg)\\
   &-\sum_{b=1}^{s-1}\bigg[\zeta_N(b;\sigma)\zeta_N(s+1-b;1)
   +(-1)^b\sum_{k=1}^{N-1}\sum_{m=1}^{N-k}\frac{\sigma^m}{k^{s+1-b}m^b}
   +\sum_{n=N+1}^{2N}\,\sum_{x=n-N}^N\frac{\sigma^x}{x^bn^{s+1-b}}\bigg].
\end{align*}
As in the proof of Proposition~\ref{prop:BBB}, we find that as $N$
grows without bound, the expression in square brackets approaches
zero when $b$ is odd, and approaches
$2\zeta(b;\sigma)\zeta(s+1-b;1)$ when $b$ is even. To complete the
proof of equation~\eqref{BBBeq75t=1}, it suffices to show that the
expressions in parentheses each tend to zero in the limit as $N$
tends to infinity.

First, since $s\ge 2$,
\begin{multline*}
  \bigg|\zeta_N(s;\sigma)\zeta_N(1;1)-\sum_{k=1}^{N-1}\sum_{m=1}^{N-k}\frac{\sigma^m}{km^s}\bigg|
  = \bigg|\sum_{k=1}^N
  \frac1k\sum_{m=N-k+1}^N\frac{\sigma^m}{m^s}\bigg|
  \le \sum_{k=1}^N \frac1k\sum_{m=N-k+1}^N \frac{1}{m^2}\\
  = \sum_{m=1}^N \frac1{m^2}\sum_{k=N-m+1}^N\frac1k
  \le \sum_{m=1}^N \frac{1}{m^2}\cdot\frac{m}{N-m+1}
  = \frac{1}{N+1}\sum_{m=1}^N\bigg(\frac1m+\frac{1}{N-m+1}\bigg)\\
  = \frac{2}{N+1}\sum_{m=1}^N\frac1m
  \to 0
\end{multline*}
as $N\to\infty$.  Also,
\begin{multline*}
  \bigg|\sum_{k=1}^{N-1}\sum_{m=1}^{N-k}\frac{\sigma^{m+k}}{k^sm}-\zeta_N(s;\sigma)\zeta_N(1;\sigma)\bigg|
  =
  \bigg|\sum_{k=1}^N\frac{\sigma^k}{k^s}\sum_{m=N-k+1}^N\frac{\sigma^m}{m}\bigg|
  \le \sum_{k=1}^N \frac1{k^2}\sum_{m=N-k+1}^N\frac1m\\
  \le \sum_{k=1}^N\frac1{k^2}\cdot\frac{k}{N-k+1}
  = \frac{1}{N+1}\sum_{k=1}^N\bigg(\frac1k+\frac{1}{N-k+1}\bigg)
  = \frac{2}{N+1}\sum_{k=1}^N\frac1k
  \to 0
\end{multline*}
as $N\to\infty$.    Finally,
\begin{multline*}
  \bigg|\sum_{n=N+1}^{2N}\,\sum_{x=n-N}^N\frac{\sigma^x}{x^sn}
  +\sum_{n=N+1}^{2N}\,\sum_{y=n-N}^N\frac{\sigma^{y+n}}{n^sy}\bigg|
  \le
  \sum_{n=N+1}^{2N}\,\sum_{x=n-N}^N\frac1{x^2n}+\sum_{n=N+1}^{2N}\,\sum_{y=n-N}^N\frac{1}{yn^2}\\
  =\sum_{x=1}^N\frac{1}{x^2}\sum_{n=N+1}^{N+x}\frac1n+\sum_{y=1}^N\frac1y\sum_{n=N+1}^{N+y}\frac1{n^2}
  \le\sum_{x=1}^N\frac1{x^2}\cdot\frac{x}{N+1}+\sum_{y=1}^N\frac1y\sum_{n=N+1}^\infty
  \frac{1}{n(n-1)}\\
  = \frac{1}{N+1}\sum_{x=1}^N\frac1x+\frac1N\sum_{y=1}^N\frac1y\to0
\end{multline*}
as $N\to\infty$.


\begin{thebibliography}{9} \raggedright
%\bibitem{AndAskRoy} G.~E.~Andrews, R.~Askey, and R.~Roy,
%\textit{Special Functions}, Cambridge University Press, 1999.

\bibitem{Berndt1} B.~Berndt, \textit{Ramanujan's Notebooks Part
I}, Springer, New York, 1985.

\bibitem{BBB} J.~M.~Borwein, D.~J.~Broadhurst, and D.~M.~Bradley,
{Evaluations of $k$-fold Euler/Zagier sums: a compendium of
results for arbitrary $k$},
%\textit{Electronic J.~Combinatorics},
\textit{Electron.~J.~Combin.}, \textbf{4} (1997), no.~2, \#R5.
Wilf Festschrift.

%\bibitem{BBG} D.~Borwein, J.~M.~Borwein, and R.~Girgensohn, Explicit
%evaluation of Euler sums, \textit{Proc.\ Edinburgh Math.\ Soc.},
%\textbf{38} (1995), 277--294.

%\bibitem{BBBLa} J.~M.~Borwein, D.~J.~Broadhurst, D.~M.~Bradley,
%and P.~Lison\v ek, Special values of multiple polylogarithms,
%\textit{Trans.\ Amer.\ Math.\ Soc.},  \textbf{353} (2001), no.~3,
%907--941. http://arXiv.org/abs/math.CA/9910045

%\bibitem{BBBLc} \bysame,
%J.~M.~Borwein, D.~J.~Broadhurst, D.~M.~Bradley, and P.~Lison\v ek,
%{Combinatorial aspects of multiple zeta values},
%\textit{Electron.~J.~Combin.}, \textbf{5} (1998), no.~1, \#R38.
%http://arXiv.org/abs/math.NT/9812020

%\bibitem{BowBrad1} D.~Bowman and D.~M.~Bradley, Resolution of some
%open problems concerning multiple zeta evaluations of arbitrary
%depth, \textit{Compositio Math.}, \textbf{139} (2003), no.~1,
%85--100. http://arXiv.org/abs/math.CA/0310061

%\bibitem{BowBradSurvey} \bysame, Multiple
%polylogarithms: A brief survey, \textit{Proceedings of a
%Conference on $q$-Series with Applications to Combinatorics,
%Number Theory and Physics}, (B.~C.~Berndt and K.~Ono eds.) Amer.\
%Math.\ Soc., Contemporary Math., \textbf{291} (2001), 71--92.
%http://arXiv.org/abs/math.CA/0310062

%\bibitem{BowBrad3} \bysame, The algebra and
%combinatorics of shuffles and multiple zeta values, \textit{J.\
%Combin.~Theory, Ser.\ A}, \textbf{97} (2002), no.~1, 43--61.
%http://arXiv.org/abs/math.CO/0310082

%\bibitem{BowBradRyoo} D.~Bowman, D.~M.~Bradley, and J.~Ryoo,
%{Some multi-set inclusions associated with shuffle convolutions
%and multiple zeta values}, \textit{European J.~Combin.},
%\textbf{24} (2003), 121--127.

%\bibitem{Prtn} D.~M.~Bradley,
%{Partition identities for the multiple zeta function}, to appear
%in \textit{Zeta Functions, Topology, and Physics}, Kinki
%University Mathematics Seminar Series, Developments in
%Mathematics. http://arXiv.org/abs/math.CO/0402091

%\bibitem{DBqMzv} \bysame,
%{Multiple $q$-zeta values}, \textit{J.\ Algebra}, \textbf{283}
%(2005), no.~2, 752--798. doi: 10.1016/j.jalgebra.2004.09.017
%http://arXiv.org/abs/math.QA/0402093



\bibitem{DBqDecomp} \bysame, {A $q$-analog of Euler's decomposition
formula for the double zeta function}, \textit{Internat.\ J.\ Math.\
Math.\ Sci.}, \textbf{2005} (2005), no.~21,  3453--3458.
doi:10.1155/IJMMS.2005.3453 [MR 2206867] (2006k:11174) {\tt
http://arxiv.org/abs/math.NT/0502002}

%\bibitem{BK1}
%D.~J.~Broadhurst and D.~Kreimer, Association of multiple zeta
%values with positive knots via Feynman diagrams up to 9 loops,
%\textit{Phys.\ Lett. B}, \textbf{393} (1997) no.~3-4, 403--412.


\bibitem{LE} L.~Euler, \textit{Meditationes Circa Singulare
Serierum Genus}, Novi Comm.\ Acad.\ Sci.\ Petropol., \textbf{20}
(1775), 140--186. Reprinted in ``Opera Omnia,'' ser.~I,
\textbf{15}, B.\ G.\ Teubner, Berlin (1927), 217--267.

\bibitem{LE2} \bysame, \textit{Briefwechsel}, vol.~1, Birh\"auser,
Basel, 1975.

\bibitem{FlajSalv} P.\ Flajolet and B.\ Salvy, Euler sums and
contour integral representations, \textit{Experiment.\ Math.},
\textbf{7} (1998), no.~1, 15--35. [MR 1618286] (99c:11110)


\bibitem{Goldbach} L.~Euler and C.~Goldbach, \textit{Briefwechsel}
1729--1764, Akademie-Verlag, Berlin, 1965.



\bibitem{Niels} N.~Nielsen, \textit{Die Gammafunktion}, Chelsea,
New York, 1965, 47--59.

\end{thebibliography}
\end{document}